\documentclass[11pt, oneside]{amsart}   	
\usepackage{geometry}                		
\usepackage{graphicx}				
\usepackage{amssymb}
\usepackage{amsmath, amscd}
\usepackage{tikz-cd}
\usepackage{tikz}

\title{On the Existence of an Orthogonal Factorization System on 1-Cob and 2-Cob}
\author{Joseph Abadi}

\begin{document}
\maketitle
\section*{Abstract}
We define the category $\bf{2-Cob}$ combinatorially and use this definition to prove the existence of an orthogonal factorization system. In the second half of the paper, we define oriented $\bf{1-Cob}$ similarly and define a functor from oriented $\bf{1-Cob}$ to $\bf{2-Cob}$. After defining this functor, the orthogonal factorization system on $\bf{2-Cob}$ is used, in turn, to prove the existence of an orthogonal factorization system on oriented $\bf{1-Cob}$. 

\section{Introduction}
An orthogonal factorization system on a category $\mathcal{C}$ is defined by two classes of morphisms $\mathcal{L}, \mathcal{R}$ such that 
\begin{enumerate}
\item The classes $\mathcal{L}$ and $\mathcal{R}$ are closed under composition and contain all isomorphisms; 
\item Any morphism $u$ factors as $me$ with $e \in \mathcal{L}, m \in \mathcal{R}$
\item Given any diagram of the following type
\begin{center}
\begin{tikzcd}
X \arrow{r}{u} \arrow{d}{e} & X' \arrow{d}{m} \\
Y \arrow{r}{v} & Y'
\end{tikzcd}
\end{center}
where $u, v$ are arbitrary morphisms and $e \in \mathcal{L}, m \in \mathcal{R}$, there exists a unique morphism $w: Y \rightarrow X'$ making the diagram commute \cite{riehl}.
\end{enumerate}
Here we find an orthogonal factorization system on the category $\bf{2-Cob}$, in which the objects are 1-manifolds and the morphisms are 2-manifolds connecting 1-manifolds, where the composition is defined by gluing of 2-manifolds. The category $\bf{2-Cob}$ is of particular interest because 2-cobordisms correspond to Feynman diagrams in string theoryÐ each boundary circle of a diagram corresponds to a closed string, and the connecting 2-manifold represents the propagation of the string from one state to another. To make use of these diagrams, it is necessary to use quantum field theory to compute the associated scattering amplitudes.

In quantum field theory, the primary objects of interest are physical fields that obey the equations of motion derived from a Lagrangian. Quantization of the field via the introduction of commutation relations causes particles to emerge as the excitations of the field. The particle states form a Hilbert space, which in particular facilitates the computation of transition amplitudes, which in turn yields probabilities for the system to transition from one state to another. Topological quantum field theories (TQFTs) are theories in which the Lagrangian depends only on the topology of a particle's worldsheet (i.e. the path it traverses in spacetime), meaning that transition amplitudes are homotopy-invariant. The canonical example of a TQFT is that given by the Chern-Simons action,
\begin{center}
\begin{equation*}
S[A] = \int\limits_{S^3} d^3x \epsilon^{\mu\nu\lambda}\textbf{Tr}(A_\mu \partial_\nu A_\lambda + A_\mu A_\nu A_\lambda)
\end{equation*}
\end{center}
where $A_\mu$ is a field in the fundamental representation of $SU(2)$ and $\epsilon^{\mu\nu\lambda}$ is the Levi-Civita symbol \cite{cs}. The action is manifestly topologically invariant because it is independent of the spacetime metric. In this case, the action describes a theory of point particles, but as we will see later, this notion can be generalized to particles in higher dimensions, such as strings.

One can use the action to explicitly calculate transition amplitudes via the path integral quantization, originally proposed by Feynman. The path integral, which lacks a formal mathematical definition, is interpreted as an integral over all possible field configurations $A(x)$. Expectation values of observables (i.e. functions of $A(x)$) are calculated by weighting the value of the function in a given field configuration by a phase given by the action associated to that configuration. For the action given above, the probability for a particle with polarization $\alpha$ at $x_1$ to propagate to a particle with polarization $\beta$ at $x_2$ is written as
\begin{center}
\begin{equation*}
\langle \Omega | A_\alpha(x_1) A_\beta(x_2) | \Omega \rangle =  \int \mathcal{D}A A_\alpha(x_1)A_\beta(x_2) \exp(iS[A])
\end{equation*}
\end{center}
In the above notation, the state $| \Omega \rangle$ denotes the vacuum of the theory. It then becomes apparent that the action serves to link physical processes to an underlying Hilbert space. Given that cobordisms correspond to Feynman diagrams of TQFTs, it is therefore natural to define a TQFT as a functor $\Phi: \bf{2-Cob} \rightarrow \bf{Hilb}$, since (as noted previously) TQFTs serve the function of assigning transition amplitudes to cobordisms via the underlying inner product structure of the state space. Roughly speaking, the orthogonal factorization system we find has as its distinguished classes of morphisms the set of processes with a single outgoing string and the set of morphisms with a single incoming string. Our orthogonal factorization system provides a systematic way to find the factorization of any cobordism and then reduces the computation of the relevant scattering amplitude to a computation of an amplitude with one incoming state and the computation of an amplitude with one outgoing state. This reduction is beneficial because in general, computations of scattering amplitudes are much simpler for diagrams in which either the initial or final state consists of a single particle \cite{qft}. 

The paper proceeds as follows. In section 2, we present the definition of $\bf{2-Cob}$ in terms of manifolds and give a heuristic description of the properties of cobordisms between $1$-dimensional manifolds. In section 3 we introduce the basic definitions of $\bf{2-Cob}$ under our interpretation (which is a combinatorial approach that we believe coincides with the traditional definition). We also define the composition of cobordisms under our definition, show it is associative, and define the cofibration and fibration classes used in our factorization system. In section 4 we prove several key lemmas needed for the proof of the existence of an orthogonal factorization system and then proceed with the theorem itself. In section 5 we give a heuristic description of oriented $\bf{1-Cob}$ (henceforth $\bf{1-Cob^+}$) and show that the factorization system on $\bf{2-Cob}$ can also be used to define an orthogonal factorization system on $\bf{1-Cob^+}$.

\section{A Heuristic Description of $\bf{2-Cob}$}
A 2-cobordism between two compact one-dimensional manifolds is defined as a tuple $(W, M, N, f, g)$, where $M$ and $N$ are (compact) 1-manifolds, $W$ is an orientable 2-manifold, and $f: M \rightarrow \partial W$, $g: N \rightarrow \partial W$ are maps such that the boundary $\partial W$ can be written as $Im(f) \sqcup Im(g)$. Intuitively, this means that a cobordism from $M$ to $N$ is specified by a 2-manifold $W$ whose boundary is $M \sqcup N$ (\cite{bartlett}). In what follows we will concern ourselves only with topological properties of the manifolds in question; that is, only the homotopy type of each manifold will be important. Due to the classification theorems for one- and two-dimensional manifolds, we can write any compact manifold $M$ as a disjoint union of circles and we can uniquely characterize any 2-manifold by its genus and number of boundaries \cite{kock}. 

With this definition in mind, it becomes evident that any connected component of a 2-cobordism can be seen as a set of ``incoming" boundary circles $M$ and a set of ``outgoing" circles $N$ that are connected by a 2-manifold $W$ that is uniquely determined by its genus. A stereotypical example of a connected component of a 2-cobordism would be the following: 
\begin{center}
\includegraphics{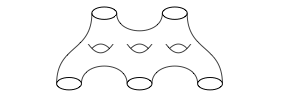}
\end{center}
In the formal definition of cobordisms we will represent a cobordism by two maps from $M$ and $N$ to a set of components, each of which is assigned a genus, so this component would be formed by mapping two elements of $M$ and three elements of $N$ to a component of genus 3. 

Cobordisms form a category with composition given by ``gluing" boundary circles. For example, the preceding example of a connected component of a cobordism could be composed with another:
\begin{center}
\includegraphics{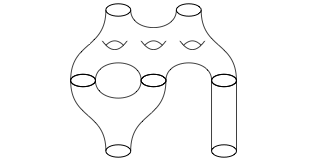}
\end{center}
It should be clear that this glued manifold is homotopic to
\begin{center}
\includegraphics{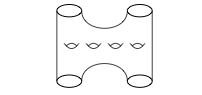}
\end{center}
which is just a connected component of a cobordism with two incoming and outgoing boundary circles and a genus of 4. This example shows that in a sense, we can think of connected components in the composition as representing equivalence classes of connected components in the original cobordism, where two individual components are considered equivalent if they share a boundary circle. Formally this will amount to defining the set of components of the composition as being a pushout of the maps from the intermediate manifold to the sets of components of the first and second cobordisms. 

One possibility that up to this point has not been discussed is the case where the manifolds $M$ and $N$ that bound $W$ are empty. In this case $W$ must of course be a 2-manifold without boundary, which we will refer to as a closed surface. Closed surfaces can, of course, be included as connected components of cobordisms between non-empty manifolds. Furthermore, in a composition of cobordisms a closed surface can be formed by connected components of the cobordisms that are not closed surfaces. One example of this is
\begin{center}
\includegraphics{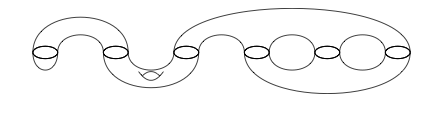}
\end{center}
which is homotopic to the closed surface of genus 3. 

This example illustrates an important point pertaining to the computation of the genus of a connected component of a composition of cobordisms. Intuitively it is clear that the boundary circles of the intermediate manifold serve to connect components of the first and second cobordisms. However, if two components of these cobordisms are already connected by a boundary circle, connecting them by a second boundary circle will form a ``hole" between the original connecting circle and the new one, and in fact every additional connecting boundary circle will form a new hole. In the example above, the component with four boundary circles is connected to the component with three boundary circles by three circles, so two holes are formed. These two holes are added to the hole included in one of the components to give a genus of three. In general, for a set of $t$ components connected by $s$ boundary circles, we add $s - t + 1$ to the genus, which we formalize in our definition of $\bf{2-Cob}$. 

The left class of morphisms in our orthogonal factorization system will consist of closed surfaces, ``merging" diagrams (which have multiple incoming components and one outgoing component), and ``annihilation" diagrams (having multiple incoming components and no outgoing components). The right class will consist of ``splitting" diagrams (which have one incoming component and multiple outgoing components) of genus zero and ``creation" diagrams (which have no incoming components and multiple outgoing components). The reason that this factorization is a somewhat natural choice (and the reason that it satisfies the properties required of an orthogonal factorization system) is that, as we will show, it greatly simplifies the calculations of pushouts and the genus in a composition of cobordisms.

\section{Elementary Properties of Cobordisms}
In this section we provide a formal definition of $\bf{2-Cob}$ in terms of finite sets mapping to a set of components, each of which is assigned a genus. Although we will not formally prove that this definition is equivalent to the definition of $\bf{2-Cob}$ in terms of manifolds, it should be clear in light of the above discussion that this combinatorial approach is actually the same.

\newtheorem*{Def1}{Definition 3.1}
\begin{Def1}
Let $X = \coprod\limits_{i \in m} S^1$, $Y = \coprod\limits_{j \in n} S^1$, where $m = \{1, \dots, m\}, n=\{1, \dots, n\}$. A $\bf{cobordism}$ $\phi: X \rightarrow Y$ is defined by specifying a finite set $C$, a function $g: C \rightarrow \mathbb{N}$, and functions $l_m: m \rightarrow C, l_n: n \rightarrow C$. We write $\phi = \{ m, n, C, l_m, l_n, g\}$. 
\begin{center}
\begin{tikzcd}
m \arrow{rd}{l_m} && n \arrow{ld}[swap]{l_n} \\
&C \arrow{d}{g} \\
& \mathbb{N}
\end{tikzcd}
\end{center}
\end{Def1}

In this definition, a connected component with $i$ incoming boundary circles, $j$ outgoing boundary circles, and genus $k$ is represented by an element $c \in C$ such that (1) $i$ elements of $m$ and $j$ elements of $n$ map to $c$, and (2) $g(c) = k$. Equipped with this interpretation, we can precisely describe several types of cobordisms. 

\newtheorem*{Def7}{Definition 3.2}
\begin{Def7}
Given a cobordism $\phi = \{m, n, C, l_m, l_n, g\}$, an element of $C$ that is not in the image of $l_m$ or $l_n$ is called a $\bf{closed}$ $\bf{surface}$ of $\phi$.
\end{Def7}

\newtheorem*{Def2}{Definition 3.3}
\begin{Def2}
A $\bf{cofibration}$ is a cobordism in which $l_n$ is injective and there exists a lift $l_n = ul_m$. 

\begin{center}
\begin{tikzcd}
m \arrow{rd}{l_m} \arrow[dashleftarrow]{rr}{\exists u} && n \arrow[hook]{ld}[swap]{l_n} \\
&C \arrow{d}{g} \\
& \mathbb{N}
\end{tikzcd}
\end{center}

A $\bf{fibration}$ is a cobordism in which $l_m$ is surjective, $l_n$ is injective, and $gl_m = 0$.

\begin{center}
\begin{tikzcd}
m \arrow[hook]{rd}{l_m} && n \arrow[two heads]{ld}[swap]{l_n} \\
&C \arrow{d}{g} \\
& \mathbb{N}
\end{tikzcd}
\end{center}
\end{Def2}

The only part of this definition which may not be immediately transparent is the requirement for the existence of a lift $u: n \rightarrow m$ in a cofibration. This condition essentially requires that any cofibration have no connected components that connect only to boundary circles of $n$. Now we move on from single cobordisms to the composition of cobordisms. The composition of cobordisms will be defined using a pushout, since this will impose an equivalence relation on components that are connected by a boundary circle of the intermediate manifold. Furthermore, as mentioned in section 2, a component with $s$ intermediate boundary circles connecting $t$ components of the original cobordisms will have a genus that is $s - t + 1$ greater than the sum of the genera of the $t$ components, since each superfluous connecting boundary circle adds a hole in the manifold.

\newtheorem*{Def4}{Definition 3.4}
\begin{Def4}
Given two cobordisms $\phi: X \rightarrow Y$, $\psi: Y \rightarrow Z$, $\phi = \{m, k, A, l_n, l_{k_1}, g_1\}$, $\psi = \{k, n, B, l_{k_2}, l_n, g_2\}$, their composition is defined by the following commutative diagram:
\begin{center}
\begin{tikzcd}
m \arrow{rd}{l_m} && k \arrow{ld}[swap]{l_{k_1}} \arrow{rd}{l_{k_2}} && n \arrow{ld}[swap]{l_n} \\
& A \arrow{rd}{\iota_A} \arrow{d}{g_1}&& B \arrow{ld}[swap]{\iota_B} \arrow{d}{g_2} \\
& \mathbb{N} & C \arrow{d}{g} & \mathbb{N} \\
&& \mathbb{N}
\end{tikzcd}
\end{center} 
where the maps $\iota_A, \iota_B$ are just those given by the pushout of $l_{k_1}, l_{k_2}$. For $x \in C$, the function $g$ is defined as 

\begin{center}
\begin{equation*}
g(x) = 1 \hspace{8 pt} + \sum\limits_{\substack{a \in A \\ \iota_A(a) = x}} (g_1(a) - 1) \hspace{8 pt} + \sum\limits_{\substack{b \in B \\ \iota_B(b) = x}} (g_2(b) - 1) \hspace{4 pt} + \sum\limits_{\substack{y \in k \\ y \mapsto x}} 1
\end{equation*}
\end{center}

\end{Def4}
We must now show that the composition is associative in order to prove that we have truly defined a category.

\newtheorem*{prop6_2}{Proposition 3.5}
\begin{prop6_2}
The composition of cobordisms is associative.
\end{prop6_2}
\begin{proof}
Suppose $\phi: W \rightarrow X, \psi: X \rightarrow Y, \chi: Y \rightarrow Z$ are cobordisms, with $\phi = \{m, n, A, a_1, \\ a_2, g_1\}, \psi = \{n, p, B, b_1, b_2, g_2\}, \chi = \{p, q, C, c_1, c_2, g_3\}$. We will denote the two different orders of composition by $\eta = \chi(\psi\phi) = \{m, q, F, h_1, h_2, g\}$ and $\eta' = (\chi\psi)\phi = \{m, q, F', h_1', h_2', g\}$. 

To determine $\eta$, we first obtain $F$ by taking the pushout $D = A \cup_n B$ and using the induced map $p \rightarrow B \rightarrow D$ to construct the pushout $F = D \cup_p C$. Similarly, for $\eta'$ we can define the pushout $E = B \cup_p C$, so we have $F' = A \cup_n E$. Since pushouts are associative (see \cite{maclane}), $F = F'$. Furthermore, this implies that the maps $m \rightarrow A \rightarrow F'$ and $m \rightarrow D \rightarrow F$ are equal; similarly, the maps $q \rightarrow C \rightarrow F$ and $q \rightarrow E \rightarrow F'$ are equal. We can then write $h_1 = h_1', h_2 = h_2'$ and draw the diagram below in an unambiguous way. 
 
\begin{center}
\begin{tikzcd}
m \arrow{dr}{a_1} && n \arrow{dl}[swap]{a_2} \arrow{dr}{b_1} && p \arrow{dl}[swap]{b_2} \arrow{dr}{c_1} && q \arrow{dl}[swap]{c_2} \\
& A \arrow{dr}{\iota_A} \arrow{d}{g_1} && B \arrow{dl}[swap]{\iota_{B_1}} \arrow{d}{g_2} \arrow{dr}{\iota_{B_2}} && C \arrow{dl}[swap]{\iota_C} \arrow{d}{g_3} \\
& \mathbb{N} & D \arrow{dr}{\iota_D} \arrow{d}{g_D} & \mathbb{N} & E \arrow{dl}[swap]{\iota_E} \arrow{d}{g_E}&\mathbb{N} \\
&&\mathbb{N}&F \arrow{d}{g,g'} &\mathbb{N} \\
&&&\mathbb{N}
\end{tikzcd}
\end{center}

It remains to check that $g(x) = g'(x)$. That is, we must check that for each element of $F$, the genus obtained via the composition $\eta$ is the same as that obtained via the composition $\eta'$.

To do so, we refer back to the formula for the genus from Definition 3.4 and expand the formula for $g_D$:

\begin{center}
\begin{equation*}
g(x) = 1 \hspace{8 pt} + \sum\limits_{\substack{d \in D \\ \iota_D(d) = x}} (g_D(d) - 1) \hspace{8 pt} + \sum\limits_{\substack{c \in C \\ \iota_E\iota_C(c) = x}} (g_3(c) - 1) \hspace{4pt} + \sum\limits_{\substack{z \in p \\ y \mapsto x}} 1 
\end{equation*}
\end{center}
We can rewrite the second term as
\begin{center}
\begin{align*}
\sum\limits_{\substack{d \in D \\ \iota_D(d) = x}} (g_D(d) - 1) & = \sum\limits_{\substack{d \in D \\ \iota_D(d) = x}} \bigg(\sum\limits_{\substack{a \in A \\ \iota_A(a) = d}} (g_1(a) - 1) \hspace{8 pt} + \sum\limits_{\substack {b \in B \\ \iota_{B_1}(b) = d}} (g_2(b) - 1) \hspace{4 pt} + \sum\limits_{\substack {y \in n \\ y \mapsto d}} 1 \bigg) \\
& =  \sum\limits_{\substack{a \in A \\ \iota_D \iota_A (a) = x}} (g_1(x) - 1) + \sum\limits_{\substack {b \in B \\ b \mapsto x}} (g_2(x) - 1) + \sum\limits_{\substack {y \in n \\ y \mapsto x}} 1
\end{align*}
\end{center} 
Similarly, we have
\begin{center}
\begin{align*}
\sum\limits_{\substack {e \in E \\ \iota_E(e) = x}} (g_E(e) - 1) &= \sum\limits_{\substack{c \in C \\ \iota_E \iota_C(c) = x}} (g_3(x) - 1) +  \sum\limits_{\substack {b \in B \\ b \mapsto x}} (g_2(x) - 1) + \sum\limits_{\substack {z \in p \\ z \mapsto x}} 1
\end{align*}
\end{center}
By using the above results, it is then easy to check that the functions
\begin{center}
\begin{equation*}
g(x) = 1 \hspace{8 pt} + \sum\limits_{\substack{d \in D \\ \iota_D(d) = x}} (g_D(d) - 1) \hspace{8 pt} + \sum\limits_{\substack{c \in C \\ \iota_E\iota_C(c) = x}} (g_3(c) - 1) \hspace{4 pt}+ \sum\limits_{\substack{z \in p \\ z \mapsto x}} 1 
\end{equation*}
\end{center}
\begin{center}
\begin{equation*}
g'(x) = 1 \hspace{8 pt} + \sum\limits_{\substack {a \in A \\ \iota_D \iota_A(a) = x}} (g_1(a) - 1) \hspace{8 pt} + \sum\limits_{\substack {e \in E \\ \iota_E(e) = x}} (g_E(e) - 1) \hspace{4 pt}+ \sum\limits_{\substack {y \in n \\ y \mapsto x}} 1
\end{equation*}
\end{center}
are equal, so in fact the formula for the genus is associative. 

\end{proof}

It will be crucial in our later analysis to know that the classes of fibrations and cofibrations are closed under composition, since in an orthogonal factorization system the left and right classes must be closed. We prove this proposition before moving on to the proof that any cobordism can be factored as a cofibration followed by a fibration.

\newtheorem*{Prop5}{Proposition 3.6}
\begin{Prop5}
The class of fibrations is closed under composition. 
\end{Prop5}
\begin{proof}
Consider the composition of two fibrations:
\begin{center}
\begin{tikzcd}
m \arrow[hook]{rd}{a} && k \arrow[two heads]{ld}[swap]{b} \arrow[hook]{rd}{c}  && n \arrow[two heads]{ld}[swap]{d} \\
& A \arrow[hook]{rd}{\iota_A} \arrow{d}{g_1}&& B \arrow[two heads]{ld}[swap]{\iota_B} \arrow{d}{g_2} \\
& \mathbb{N} & C \arrow{d}{g} & \mathbb{N} \\
&& \mathbb{N}
\end{tikzcd}
\end{center} 
where $g_1a = 0$ and $g_2c = 0$. The map $\iota_A$ is an injection because $c$ is injective, and $\iota_B$ is a surjection because $b$ is surjective. Then $\iota_A a$ and $\iota_B d$ are injective and surjective, respectively, as desired. We must also show that $g \iota_A a = 0$ and that there are no closed surfaces in $C$. There are no closed surfaces because $\iota_B d$ is surjective. To see that $g \iota_A a = 0$, consider the formula for $g$:
\begin{center}
\begin{equation*}
g(x) = 1 \hspace{8 pt} + \sum\limits_{\substack{a \in A \\ \iota_A(a) = x}} (g_1(a) - 1) \hspace{8 pt} + \sum\limits_{\substack{b \in B \\ \iota_B(b) = x}} (g_2(b) - 1) \hspace{4 pt} + \sum\limits_{\substack{y \in k \\ y \mapsto x}} 1
\end{equation*}
\end{center}
We assume $x = \iota_A a(y)$ for some $y \in m$. If we also have $x \in im(\iota_B)$, then the $g_2$ terms are all zero, as $x = \iota_B c(z)$ for some $z \in k$ and $g_2 c = 0$. Likewise, the $g_1$ terms are zero. We are left with 
\begin{center}
\begin{equation*}
g(x) = \sum\limits_{\substack{b \in B \\ \iota_B(b) = x}} (-1) \hspace{4 pt} + \sum\limits_{\substack{y \in k \\ y \mapsto x}} 1
\end{equation*}
\end{center}
since $\iota_A$ is injective. But this sum is clearly zero, since every element of $\iota_B$ that maps to $x$ must have come from $k$ given that $x \in im(\iota_A)$, so $g\iota_A a = 0$ in this case. 

When $x \notin im(\iota_B)$, we have 
\begin{center}
\begin{equation*}
g(x) = 1 \hspace{8 pt} + \sum\limits_{\substack{a \in A \\ \iota_A(a) = x}} (g_1(a) - 1) 
\end{equation*}
\end{center}
which is zero because there is only one element in $\iota_A$ mapping to $x$, and the genus of that element is zero because it is in $im(a)$. Therefore the composition of two fibrations is again a fibration.  
\end{proof}

Note that Proposition 5 also shows that the class of cofibrations is closed under composition, since the proof is just a special case of the proof for fibrations. 

\section{Factorization of Cobordisms}
Now we are nearly ready to show that any cobordism factors uniquely as a cofibration followed by a fibration. Before doing so we prove some important lemmas.

\newtheorem*{Lemma1}{Lemma 4.1}
\begin{Lemma1}
A commutative square in which all maps are injective is a pushout square if and only if it is a pullback square. 
\end{Lemma1}
\begin{center}
\begin{tikzcd}
A \arrow[hook]{r}{u} \arrow[hook]{d}{v} & B \arrow[hook]{d}{s}\\
C \arrow[hook]{r}{t}& D
\end{tikzcd}
\end{center}
\begin{proof}
For one direction of the proof, we must show that $A$ is isomorphic to $B \times_D C$ when the diagram above is a pushout square. Consider the function $\phi: A \rightarrow B \times_D C$ given by $\phi(x) = (u(x), v(x))$ (which is well-defined because the square is commutative). If $(y, z) \in B \times_D C$, since $D = B \sqcup C / \sim$, where $a \sim b$ if there exists $c \in A$ such that $u(c) = a, v(c) = b$, it must be that there is $w \in A$ such that $u(w) = y, v(w) = z$. This $w$ is unique because $u, v$ are injections. Therefore $\phi$ is injective and surjective, so there is an isomorphism between $\phi: A \rightarrow B \times_D C$ given by $\phi(x) = (u(x), v(x))$, meaning $A$ is in fact the pullback.

The other direction of the proof is analogous.

\end{proof}

\newtheorem*{Lemma2}{Lemma 4.2}
\begin{Lemma2}
Consider the following composition of cobordisms:
\begin{center}
\begin{tikzcd}
m \arrow{rd}{l_m} \arrow[dashleftarrow]{rr}{\exists u} && k \arrow[hook]{ld}[swap]{l_{k}} \arrow{rd}{j_k} && n \arrow{ld}[swap]{j_n} \\
& A \arrow{rd}{\iota_A} && B \arrow[hook]{ld}[swap]{\iota_B} \\
&& C
\end{tikzcd}
\end{center}
In this case the set of closed surfaces of $C$, Cl$(C)$, is in bijection with Cl$(A) \sqcup $Cl$(B)$.
\end{Lemma2}
\begin{proof}
The closed surfaces of $C$ are precisely those elements of $C$ that are not in the image of $\iota_A l_m$ or $\iota_B j_n$. On the other hand, the closed surfaces of $A$ and $B$ are elements of $A$ and $B$ not in the images of $l_m, l_k$ or $j_k, j_n$, respectively. 

Suppose $x \in C$. Then, by the definition of the pushout, there is some $y \in A$ or $z \in B$ mapping to $x$. If there is $y$ such that $\iota_A(y) = x$ and $y$ is not a closed surface, there must be some $a \in m$ or $b \in k$ mapping to $y$. In the former case $x$ is not closed, and in the latter, the lift $l_mu(b) = y$ shows that $x$ is in the image of $\iota_A l_m$ regardless. Therefore any such $y$ must be a closed surface. Now suppose there is $z$ such that $\iota_B(z) = x$. If $z$ is not closed, it is either in the image of $j_n$, in which case $x$ is not closed, or there is some $b \in k$ such that $j_k(b) = z$. But then $\iota_Al_m(u(b)) = x$, so $x$ is not a closed surface. Therefore Cl$(A) \sqcup $Cl$(B)$ surjects onto Cl$(C)$. 

On the other hand, suppose there are two closed surfaces $y \in A$, $z \in B$ mapping to the same $x \in C$. By definition, there must be some $b_1, b_2 \in k$ such that $l_k(b) = y, j_k(b) = z$. But this is impossible because $y, z$ are closed surfaces. The map Cl$(A) \sqcup $Cl$(B) \rightarrow $Cl$(C)$ given by $(x, i) \rightarrow \iota_i(x)$ is therefore an injection as well, so it is an isomorphism, as desired. 
\end{proof}

Now we are prepared to prove one of the main propositions needed to demonstrate the existence of an orthogonal factorization system on $\bf{2-Cob}$.

\newtheorem*{Prop3}{Proposition 4.3}
\begin{Prop3}
Every cobordism factors uniquely as a cofibration followed by a fibration.
\end{Prop3}
\begin{proof}
Consider a cobordism
\begin{center}
\begin{tikzcd}
m \arrow{dr}{a} && n \arrow{dl}[swap]{d} \\
&C \arrow{d}{g} \\
&\mathbb{N}
\end{tikzcd}
\end{center}
It can be decomposed as
\begin{center}
\begin{tikzcd}
m \arrow{rd}{a} \arrow[dashleftarrow]{rr}{\exists u} && k \arrow[hook]{ld}[swap]{b} \arrow[hook]{rd}{c} && n \arrow[two heads]{ld}[swap]{d} \\
& im(a) \sqcup (C - im(b) - im(a)) \arrow[hook]{rd}{\iota_1} \arrow{d}{g_1} && im(b) \arrow[hook]{ld}[swap]{\iota_2} \arrow{d}{g_2} \\
&\mathbb{N} & C \arrow{d}{g} & \mathbb{N} \\
&& \mathbb{N}
\end{tikzcd}
\end{center}
where $a, d, \iota_1, \iota_2$ are the canonical maps given by the original diagram, $k$ is the pullback of the injections $\iota_1, \iota_2$ (which then forms a pushout square by Lemma 1), and $g_1, g_2$ agree with the restrictions of $g$ to the sets $im(a) \sqcup (C - im(b) - im(a))$, $im(b)$, respectively. By inspection, it is easy to check that $g_1$ and $g_2$ define the same $g$ as in the diagram above.

As for uniqueness, we define two factorizations 
\begin{center}
\begin{tikzcd}
m \arrow{rd}{a} \arrow[dashleftarrow]{rr}{\exists u} && k \arrow[hook]{ld}[swap]{b} \arrow[hook]{rd}{c} && n \arrow[two heads]{ld}[swap]{d} \\
& A \arrow[hook]{rd}{\iota_A} \arrow{d}{g_1} && B \arrow[hook]{ld}[swap]{\iota_B} \arrow{d}{g_2} \\
&\mathbb{N} & C \arrow{d}{g} & \mathbb{N} \\
&& \mathbb{N}
\end{tikzcd}
\begin{tikzcd}
m \arrow{rd}{a'} \arrow[dashleftarrow]{rr}{\exists u'} && k' \arrow[hook]{ld}[swap]{b'} \arrow[hook]{rd}{c'} && n \arrow[two heads]{ld}[swap]{d'} \\
& A' \arrow[hook]{rd}{\iota_A'} \arrow{d}{g_1'} && B' \arrow[hook]{ld}[swap]{\iota_B'} \arrow{d}{g_2'} \\
&\mathbb{N} & C \arrow{d}{g} & \mathbb{N} \\
&& \mathbb{N}
\end{tikzcd}
\end{center}
to be equivalent if there are isomorphisms $\psi: A \simeq A'$, $\phi: B \simeq B'$, $\chi: k \simeq k'$ such that $\psi a = a'$, $\phi d = d'$, $\chi b = b' \psi$, $\chi c = c' \phi$. 

Note first that $B \simeq B'$, since $\iota_B d = \iota_B' d'$ and we have two epimorphism-monomorphism factorizations of this map, 
\begin{center}
\begin{tikzcd}
& B \arrow[hook]{rd} \\
n \arrow[two heads]{ru} \arrow[two heads]{rd} && C \\
& B' \arrow[hook]{ru}
\end{tikzcd}
\end{center}
so there is an isomorphism $\phi: B \rightarrow B'$ that can be picked such that $\phi d = d'$. 

For $x \in A'$, set $\psi(x) = \iota_{A}^{-1} \iota_{A'} (x)$. We need to show that this is well-defined. 
There are two cases: $x \in im(a')$ and $x \notin im(a')$. In the first this is clearly well-defined, since there is some $y \in m$ such that $a'(y) = x$, so given that $\iota_{A'}a' = \iota_A a$, we have $\psi(x) = \iota_A^{-1}\iota_{A'} a'(y) = \iota_A^{-1} \iota_A a(y) = a(y)$. In the second case, given that $b'$ lifts to $a'u'$, it must be that $x$ is a closed surface. Then Lemma 1 shows that Cl$(A') \sqcup $Cl$(B') \simeq C \simeq $Cl$(A) \sqcup $Cl$(B)$, but Cl$(B), $Cl$(B')$ are both empty because $d'$ is surjective, meaning that Cl$(A) \simeq $Cl$(A') \simeq C$. Therefore $\iota_{A'}(x) \in im(\iota_A)$, so again the map is well-defined. But then $\psi$ and $\psi^{-1}$ are mutually inverse on finite sets, hence are isomorphisms. It is easily checked that $\psi$ satisfies the required condition:
\begin{center}
$\psi a'(x) = \iota^{-1}_A \iota_{A'} a'(x) = \iota_A^{-1} \iota_A a(x) = a(x)$
\end{center}
Then $A \simeq A'$, $B \simeq B'$. But by Lemma 1, then, $k$ and $k'$ are pullbacks of the same pair of maps, so they are isomorphic as well by the universality of the pullback. The spans $A \leftarrow k \rightarrow B$, $A' \leftarrow k' \rightarrow B'$ have the same pushout, so we can find an isomorphism $\chi: k \rightarrow k'$ such that $\chi b = b' \psi, \chi c = c' \phi$. Therefore the factorizations are the same up to genus. 

Now consider the formula for the genus from Definition 3.4 for computing the function $g$:
\begin{center}
\begin{equation*}
g(x) = 1 \hspace{8 pt} + \sum\limits_{\substack{a \in A \\ \iota_A(a) = x}} (g_1(a) - 1) \hspace{8 pt} + \sum\limits_{\substack{b \in B \\ \iota_B(b) = x}} (g_2(b) - 1) \hspace{4 pt} + \sum\limits_{\substack{y \in k \\ y \mapsto x}} 1
\end{equation*}
\end{center}
We will show that for any $x$, only one term of either of the sums involving $g_1$ or $g_2$ can be nonzero, so it is uniquely determined by g(x). There are three cases: either $x \in im(\iota_A) \cap im(\iota_B)$, $x \in im(\iota_A) - im(\iota_B)$, or $x \notin im(\iota_A)$. 

In the first case, when $x \in im(\iota_B)$, it must be that $x \in im(\iota_B d)$ (since $d$ is surjective), in which case the element $b \in B$ mapping to $x$ must satisfy $g_2(b) = 0$, as the second cobordism is a fibration. There can only be one nonzero $g_1$ term because $\iota_A$ is an injection.

In the second case, clearly there are no $g_2$ terms, and again there can be only one nonzero $g_1$ term because $\iota_A$ is an injection.

The third case is the same as the second case, since there are now no $g_1$ terms and only one possible nonzero $g_2$ term.

Therefore the functions $g_1, g_2$ are uniquely determined by $g$ and the pullback square in the diagram, so we conclude that any two factorizations must have the same functions $g_1, g_2$. Then any two factorizations of this cobordism are equivalent. 


\end{proof}

\newtheorem*{Prop4}{Proposition 4.4}
\begin{Prop4}
Cofibrations are epimorphisms in \bf{2-Cob}.
\end{Prop4}
\begin{proof}
Given a cobordism $\chi$ and a cofibration $\psi$, suppose there are two cobordisms $\phi, \phi'$ such that $\phi\psi = \phi'\psi = \chi$. Diagramatically, this is
\begin{center}
\begin{tikzcd}
m \arrow{rd}{a} \arrow[dashleftarrow]{rr}{\exists u} && k \arrow[hook]{ld}[swap]{b} \arrow{rd}{c} && n \arrow{ld}[swap]{d} \\
& A \arrow{rd}{\iota_A} \arrow{d}{g_1}&& B \arrow[hook]{ld}[swap]{\iota_B} \arrow{d}{g_2} \\
& \mathbb{N} & C \arrow{d}{g} & \mathbb{N} \\
&& \mathbb{N}
\end{tikzcd}
\begin{tikzcd}
m \arrow{rd}{a} \arrow[dashleftarrow]{rr}{\exists u} && k \arrow[hook]{ld}[swap]{b} \arrow{rd}{c'} && n \arrow{ld}[swap]{d'} \\
& A \arrow{rd}{\iota_A}\arrow{d}{g_1} && B' \arrow[hook]{ld}[swap]{\iota_B'} \arrow{d}{g_2} \\
& \mathbb{N} & C \arrow{d}{g} & \mathbb{N} \\
&& \mathbb{N}
\end{tikzcd}
\end{center}
with $\iota_B d = \iota_B' d'$.

We must find an isomorphism $\phi: B \simeq B'$ such that $\phi d = d'$. We will set $\phi = \iota_{B'}^{-1} \iota_B$ and show that it is well-defined. There are three cases to consider for $x \in B$: $x \in im(d)$; $x \notin im(d), \iota_B(x) \in im(\iota_A)$; and $x \notin im(d), \iota_B(x) \notin im(\iota_A)$.

In the first case, we have some $y \in n$ such that $d(y) = x$, and $\iota_B d = \iota_{B'} d'$ by assumption, so $\iota_{B'}^{-1} \iota_B (x) = \iota_{B'}^{-1}\iota_{B'} d'(y)$, which is well-defined because $\iota_{B'}$ is injective. 

In the second case, by assumption there is $y \in A$ such that $\iota_A(y) = \iota_B(x)$. By the definition of the pushout, there are $z, z' \in k$ such that $b(z) = y, c(z') = x$. But then, since $\iota_{B'} c' = \iota_A b = \iota_B c$, we have $\iota_{B'}^{-1}\iota_B(x) = \iota_{B'}^{-1} \iota_{B'} c(z')$, which again is well-defined. 

In the last case, $x$ must map to a closed surface of $C$. By Lemma 2, 
\begin{center}
Cl$(C) \simeq $Cl$(A) \sqcup $Cl$(B) \simeq $Cl$(A) \sqcup $Cl$(B')$
\end{center}
but we know that $\iota_B(x) \notin im(\iota_A)$, so it must be that $\iota_B(x) \in im(\iota_{B'})$ as well. 

The proof that $\phi$ is an isomorphism is immediate by checking that $\phi$ and $\phi^{-1}$ are mutually inverse, where $\phi^{-1} = \iota_B^{-1} \iota_{B'}$. It is furthermore clear that $\phi d = d'$. 

Now we check that $g_2 = g_2'$. The formula for computing $g$ is
\begin{center}
\begin{equation*}
g(x) = 1 \hspace{8 pt} + \sum\limits_{\substack{a \in A \\ \iota_A(a) = x}} (g_1(a) - 1) \hspace{8 pt} + \sum\limits_{\substack{b \in B \\ \iota_B(b) = x}} (g_2(b) - 1) \hspace{4 pt} + \sum\limits_{\substack{y \in k \\ y \mapsto x}} 1
\end{equation*}
\end{center}
For each $x$, there can be only one nonzero term coming from the sum involving $g_2$ because $\iota_B$ is injective. This term is uniquely determined because the other terms are all known. Therefore there is only one possible function $g_2$ such that the composition of the cobordisms yields $\chi$. Then $\phi = \phi'$, so cofibrations are epimorphisms on $\bf{2-Cob}$.
\end{proof}

The proof of Proposition 4.4 also shows that fibrations are epimorphisms on $\bf{2-Cob}$, as this is just a special case of that proposition. Now we have enough to prove the main theorem:

\newtheorem*{Thm6}{Theorem 4.5}
\begin{Thm6}
There exists an orthogonal factorization system on $\bf{2-Cob}$.
\end{Thm6}
\begin{proof}
In the diagram below (where $e$ is a cofibration and $m$ is a fibration), note that $mu = ve = gf$ for a unique cofibration-fibration pair $f, g$ by Proposition 4.3. Also, $u$ and $v$ factor uniquely as $u = ba$, $v = dc$. Therefore $gf = dce = mba$, so $g = d = mb, f = ce = a$. By Proposition 4.4, $c$ is the unique cofibration such that $f = ce$ and $b$ is the unique fibration such that $g = mb$. Set $w = bc$. Then $mwe = (mb)(ce) = gf = ve = mu$, as desired. 

Now suppose there exists some $w'$ such that $mw'e = ve = mu$. Since $w'$ factors uniquely as $w' = rs$ by Proposition 4.3, $mwe = mrse = gf$, we also have $mr = g, se = f$. By Proposition 4.4 again, $r = b, s = c$, so $w'$ is uniquely determined. 
\end{proof}

\begin{center}
\begin{tikzcd}
X \arrow{r}{u} \arrow{d}{e} & X' \arrow{d}{m} \\
Y \arrow{r}{v} \arrow[dotted]{ur}{\exists !} & Y'
\end{tikzcd}
\end{center}

\section{Orthogonal Factorization on Oriented 1-Cob}
We begin in Section 5.1 by giving an intuitive description of the properties of $\bf{1-Cob^+}$ and draw parallels to the case of $\bf{2-Cob}$. In Sections 5.2, 5.3, and 5.4 we give a definition of $\bf{1-Cob^+}$ and show that there is a functor $\Phi: \bf{1-Cob^+} \rightarrow \bf{2-Cob}$. We proceed to prove that the orthogonal factorization on $\bf{2-Cob}$ yields a factorization system on $\bf{1-Cob^+}$. Our main task in proving this statement is showing that for any cobordism $\phi \in im(\Phi)$, if $\phi$ factors as a cofibration-fibration pair $\chi\psi$, then $\chi, \psi \in im(\Phi)$, but we will also prove several minor lemmas needed to complete the full proof. 

\subsection{A Heuristic Description of Oriented 1-Cob}
The general definition of an $n$-cobordism outlined previously shows that a 1-cobordism is just a pair of finite sets $M, N$ with a curve connecting each point of $M \sqcup N$ to some other point in that set and a certain number of closed loops making up the rest of the manifold $W$ connecting $M$ and $N$. An oriented 1-cobordism differs from a general 1-cobordism in that there are certain restrictions on how curves can be drawn between points of $M \sqcup N$. In $\bf{1-Cob^+}$, every point of the sets $M$ and $N$ is assigned a sign, so we can write $M = M_+ + M_-$ and $N = N_+ + N_-$. The curves that make up the 1-manifold must have one endpoint in $M_+ + N_-$ and the other in $M_- + N_+$, so this construction essentially yields a bijection $M_+ + N_- \rightarrow M_- + N_+$. An example of an oriented 1-cobordism would be
\begin{center}
\begin{tikzpicture}
\draw[fill]  (0,0) circle [radius = 0.1];
\draw[fill] (1,0) circle [radius = 0.1];
\draw[fill] (0,3) circle [radius = 0.1];
\draw[fill] (1,3) circle [radius = 0.1];
\draw[fill] (2,0) circle [radius = 0.1];
\draw[fill] (3,0) circle [radius = 0.1];
\draw[fill] (2,3) circle [radius = 0.1];
\draw[fill] (3,3) circle [radius = 0.1];
\draw[fill] (4,3) circle [radius = 0.1];
\draw[fill] (5,3) circle [radius = 0.1];
\draw[fill] (6,3) circle [radius = 0.1];
\draw[fill] (7,3) circle [radius = 0.1];
\node [below] at (0,0) {+}; 
\node [below] at (1,0) {\textendash};
\node [below] at (2,0) {+};
\node [below] at (3,0) {\textendash};
\node [above] at (0,3) {+};
\node [above] at (1,3) {+};
\node [above] at (2,3) {\textendash};
\node [above] at (3,3) {+};
\node [above] at (4,3) {\textendash};
\node [above] at (5,3) {\textendash};
\node[above] at (6,3) {+};
\node[above] at (7,3) {\textendash};
\draw[thick] (0,0) -- (0,3);
\draw[thick] (1,0) -- (5,3);
\draw[thick] (3,0) arc [radius = 0.5, start angle = 0, end angle = 180];
\draw[thick] (1,3) arc [radius = 0.5, start angle = -180, end angle = 0];
\draw[thick] (3,3) arc [radius = 0.5, start angle = -180, end angle = 0];
\draw[thick] (6,3) arc [radius = 0.5, start angle = -180, end angle = 0];
\end{tikzpicture}
\end{center}
Composition is again defined by gluing. Although it is not immediately apparent, gluing two oriented cobordisms together does indeed result in a 1-cobordism after ``straightening" lines. An example of a composition involving the above cobordism is
\begin{center}
\begin{tikzpicture}
\draw[fill]  (0,0) circle [radius = 0.1];
\draw[fill] (1,0) circle [radius = 0.1];
\draw[fill] (0,3) circle [radius = 0.1];
\draw[fill] (1,3) circle [radius = 0.1];
\draw[fill] (2,0) circle [radius = 0.1];
\draw[fill] (3,0) circle [radius = 0.1];
\draw[fill] (2,3) circle [radius = 0.1];
\draw[fill] (3,3) circle [radius = 0.1];
\draw[fill] (4,3) circle [radius = 0.1];
\draw[fill] (5,3) circle [radius = 0.1];
\draw[fill] (6,3) circle [radius = 0.1];
\draw[fill] (7,3) circle [radius = 0.1];
\draw[fill] (0,6) circle [radius = 0.1];
\draw[fill] (1,6) circle [radius = 0.1];
\draw[fill] (2,6) circle [radius = 0.1];
\draw[fill] (3,6) circle [radius = 0.1];
\node [below] at (0,0) {+}; 
\node [below] at (1,0) {\textendash};
\node [below] at (2,0) {+};
\node [below] at (3,0) {\textendash};
\node [left] at (0,3) {+};
\node [left] at (1,3) {+};
\node [left] at (2,3) {\textendash};
\node [left] at (3,3) {+};
\node [left] at (4,3) {\textendash};
\node [left] at (5,3) {\textendash};
\node[left] at (6,3) {+};
\node[left] at (7,3) {\textendash};
\node[above] at (0,6) {+};
\node[above] at (1,6) {+};
\node[above] at (2,6) {\textendash};
\node[above] at (3,6) {\textendash};
\draw[thick] (0,0) -- (0,3);
\draw[thick] (1,0) -- (5,3);
\draw[thick] (3,6) --(5,3);
\draw[thick] (0,6) -- (0,3);
\draw[thick] (1,6) -- (1,3);
\draw[thick] (2,6) -- (4,3);
\draw[thick] (3,0) arc [radius = 0.5, start angle = 0, end angle = 180];
\draw[thick] (1,3) arc [radius = 0.5, start angle = -180, end angle = 0];
\draw[thick] (3,3) arc [radius = 0.5, start angle = -180, end angle = 0];
\draw[thick] (6,3) arc [radius = 0.5, start angle = -180, end angle = 0];
\draw[thick] (3,3) arc [radius = 0.5, start angle = 0, end angle = 180];
\draw[thick] (7,3) arc [radius = 0.5, start angle = 0, end angle = 180];
\end{tikzpicture}
\end{center}
By following each line from endpoint to endpoint, it becomes apparent that this composition is the cobordism
\begin{center}
\begin{tikzpicture}
\draw[fill] (0,3) circle [radius = 0.1];
\draw[fill] (1,3) circle [radius = 0.1];
\draw[fill] (2,3) circle [radius = 0.1];
\draw[fill] (3,3) circle [radius = 0.1];
\draw[fill] (0,0) circle [radius = 0.1];
\draw[fill] (1,0) circle [radius = 0.1];
\draw[fill] (2,0) circle [radius = 0.1];
\draw[fill] (3,0) circle [radius = 0.1];
\draw[thick] (0,0) -- (0,3);
\draw[thick] (1,0) -- (3,3);
\draw[thick] (1,3) arc [radius = 0.5, start angle = -180, end angle = 0];
\draw[thick] (3,0) arc [radius = 0.5, start angle = 0, end angle = 180];
\draw[thick] (4.5,1.5) circle [radius=0.5];
\node[above] at (0,3) {+};
\node[above] at (1,3) {+};
\node[above] at (2,3) {\textendash};
\node[above] at (3,3) {\textendash};
\node[below] at (0,0) {+};
\node[below] at (1,0) {\textendash};
\node[below] at (2,0) {+};
\node[below] at (3,0) {\textendash};
\end{tikzpicture}
\end{center}
This example serves to point out two important facts about compositions of oriented 1-cobordisms. In a composition $M \xrightarrow{\phi} K \xrightarrow{\psi} N$, to understand where a point of $M_+ + N_-$ is mapped, one must first find the point it is paired with under $\phi$, then find the point matched to this new point under $\psi$, and if that point is still not in $M_- + N_+$, find the point to which that one is matched under $\phi$, and so forth until reaching a point of $M_- + N_+$. The second fact is that closed loops can be formed in a composition in which neither factor explicitly includes a closed loop, in analogy with the situation in $\bf{2-Cob}$ in which closed surfaces are formed in the composition. 

The notion of having to follow a line from endpoint to endpoint and the formation of closed loops in the composition are both captured in our definition of $\bf{1-Cob^+}$ by the concept of the pushout. Intuitively speaking, this is because following lines through the composition essentially entails establishing an equivalence relation in which equivalence classes are points that fall on the same line or loop, which is precisely the equivalence relation defined by the pushout. In fact, our definition of $\bf{1-Cob^+}$ will be analogous to our definition of $\bf{2-Cob}$ with two key differences: first, there will be no function for the genus because 1-dimensional manifolds have no property analogous to genus, and second, there will be serious restrictions on what may map to a connected component of the cobordism corresponding to the fact that the cobordism induces a bijection of signed sets. 

Just as in $\bf{2-Cob}$, the category $\bf{1-Cob^+}$ has cobordisms that can be labeled as fibrations or cofibrations. The class of fibrations is made up of cobordisms $M \rightarrow N$ in which no two points of $M$ are connected to each other (i.e. no particles in $M$ ``annihilate'') and there are no closed loops. The class of cofibrations consists of diagrams in which no points of $N$ are connected to each other (i.e. no particle-antiparticle pairs are created). When we show that the orthogonal factorization system on $\bf{2-Cob}$ implies the existence of an orthogonal factorization system on $\bf{1-Cob}$, we will see that in fact these are precisely the morphisms that map to fibrations and cofibrations in $\bf{2-Cob}$ under a functor $\\ \Phi: \bf{1-Cob^+} \rightarrow \bf{2-Cob}$ (to be defined later).
  
\subsection{Properties of Oriented 1-Cob}
\newtheorem*{def3.1}{Definition 5.2.1}
\begin{def3.1}
In $\bf{1-Cob^+}$, a cobordism of finite sets $\theta: m = m_- + m_+ \rightarrow n_- + n_+ = n$ is defined by a diagram of the following form:
\begin{center}
\begin{tikzcd}
m \arrow{dr}{l_m} && n \arrow{dl}[swap]{l_n} \\
& C
\end{tikzcd}
\end{center}
where each $c \in C$ satisfies one of the following:
\begin{itemize}
\item $l_m^{-1}(c) = \{a, b\}, a \in m_+, b \in m_-, l_n^{-1}(c) = \emptyset $
\item $l_n^{-1}(c) = \{a, b\}, a \in n_+, b \in n_-, l_m^{-1}(c) = \emptyset$ 
\item $l_m^{-1}(c) = \{a\}, l_n^{-1}(c) = \{b\}, a \in m_+, b \in n_+$ 
\item $ l_m^{-1}(c) = \{a\}, l_n^{-1}(c) = \{b\}, a \in m_-, b \in n_- $
\item $l_m^{-1}(c) = l_n^{-1}(c) = \emptyset$

\end{itemize}
\end{def3.1}

Given two cobordisms $\{m, k, A, l_m, l_k\}, \{k, n, B, l_k', l_n\}$ we define their composition by 

\begin{center}
\begin{tikzcd}
m \arrow{dr}{l_m} && k \arrow{ld}[swap]{l_k} \arrow{rd}{l_k'} && n \arrow{ld}[swap]{l_n} \\
& A \arrow{rd}{\iota_A} && B \arrow{ld}[swap]{\iota_B} \\
&&C
\end{tikzcd}
\end{center}

where $C$ is the pushout of the span $A \leftarrow k \rightarrow B$. We must show that this composition is well-defined; that is, that the composition of two oriented 1-cobordisms yields an oriented 1-cobordism. To do so, it suffices to show that this definition of $\bf{1-Cob^+}$ is equivalent to the definition in \cite{abram} of the free closed compact category on one element (which is equivalent to $\bf{1-Cob^+}$ \cite{freyd}). The formal proof of this equivalence is rather tedious, so we include only a sketch here.

In \cite{abram}, the free compact category on one element is defined as the category in which objects are signed sets $A = \{A_+, A_-\}$ and a morphism between two objects $A$ and $B$ is a bijection $A_+ + B_- \rightarrow A_- + B_+$ and a number $n$ of ``loops." The bijection $A_+ + B_- \rightarrow A_- + B_+$ is interpreted in the following way: the cardinalities of the two sets must be the same, so $| A_+ + B_- | = | A_- + B_+ | = k$ for some $k \in \mathbb{N}$. Then, by numbering each element of each set, the bijection can be written as a permutation $\pi \in S(k)$. This means that each morphism can be characterized uniquely as $(n, \pi)$ where $n$ is again the number of loops and $\pi \in S(k)$. 

With this interpretation it is easy to define the composition of two morphisms $(A \rightarrow B) = (n, \pi)$, $(B \rightarrow C) = (m, \sigma)$. The permutations $\pi$ and $\sigma$ can clearly be broken up into matrices in the following form:
\begin{center}
\[ \pi = 
\begin{pmatrix}
\pi_{A_+A_-} & \pi_{A_+B_+} \\
\pi_{B_-A_-} & \pi_{B_-B_+}
\end{pmatrix}
\hspace{20pt}
\sigma = 
\begin{pmatrix}
\sigma_{B_+B_-} & \sigma_{B_+C_+} \\
\sigma_{C_-B_-} & \sigma_{C_-C_+}
\end{pmatrix}
\]
\end{center}

The composition of permutations is given by the ``execution formula" $\theta = $Ex$(\pi, \sigma)$, 
\begin{center}
\[ Ex(\pi, \sigma) = \theta = 
\begin{pmatrix}
\theta_{A_+A_-} & \theta_{A_+C_+} \\
\theta_{C_-A_-} & \theta_{C_-C_+}
\end{pmatrix}
\]
\end{center}
where, denoting composition by a semicolon and indefinite repetition by an asterisk, the components of $\theta$ are defined by
\begin{center}
\begin{align*}
\theta_{A_+A_-} &= \pi_{A_+A_-} \cup \pi_{A_+B_+} ; \sigma_{B_+B_-}; (\pi_{B_-B_+}; \sigma_{B_+B_-})^*; \pi_{B_- A_-} \\
\theta_{A_+C_+} &= \pi_{A_+B_+} ; (\sigma_{B_+B_-}; \pi_{B_-B_+})^*; \sigma_{B_+C_+} \\
\theta_{C_-A_-} &= \sigma_{C_-B_-}; (\pi_{B_-B_+}; \sigma_{B_+B_-})^*; \pi_{B_-A_-} \\
\theta_{C_-C_+} &= \sigma_{C_-C_+} \cup \sigma_{C_-B_-}; (\pi_{B_-B_+}; \sigma_{B_+B_-})^*; \pi_{B_-B_+}; \sigma_{B_+C_+} 
\end{align*}
\end{center}

Intuitively, this formula says that $\theta$ should be obtained by first applying $\pi$ to points of $A_+$ and $\sigma$ to points of $C_-$ and then continuing to apply the appropriate permutations until reaching a point of $A_- + C_+$. In \cite{abram}, the author proves a proposition regarding the geometry of $S(n)$ that shows this formula is always well-defined. There is a similar formula to determine the loops $\mathcal{L}(\pi, \sigma)$ formed in the composition. Newly formed loops are equivalence classes in $B$, and a point $j \in B$ is in a loop if it is fixed by some permutation of the form $(\pi_{B_-B_+}; \sigma_{B_+B_-})^*$ or $(\sigma_{B_+B_-};\pi_{B_-B_+})^*$. Adding the number of newly formed loops to the number of loops in the original two morphisms yields the number of loops in the composition. 

With this understanding of the definition of $\bf{1-Cob^+}$ given in \cite{abram}, we can proceed to give a heuristic argument as to why our definition is equivalent. First note that in our definition, a cobordism from $m = m_+ + m_-$ to $n = n_+ + n_-$ does indeed induce a bijection $m_+ + n_- \rightarrow m_- + n_+$. This is because of the conditions imposed on the functions $l_m$ and $l_n$, which state that either the preimage of every point in $C$ is either empty or consists of exactly one point from both sets. The points with empty preimage correspond to the loops in the definition of \cite{abram}, so we can in fact interpret every oriented 1-cobordism under our definition as a permutation and a number of loops. 

It remains to understand why our definition of the composition is equivalent to that in \cite{abram}. Consider the composition of $\{m, k, A, l_m, l_k\}$ with $\{k, n, B, l_k', l_n\}$. The bijection $m_+ + n_- \rightarrow m_- + n_+$ is obtained by first determining the pushout $A \cup_k B$. The pushout is defined by an equivalence relation $\sim$ such that $a \sim b$ if there is $z \in k$ such that $l_k(z) = a$, $l_k'(z) = b$. Then, for example, a point $x \in m_+$ maps to the same component of $A \cup_k B$ as a point $y \in m_-$ if and only if one of the following holds: 
\begin{enumerate}
\item The two points map to the same component under $l_m$; that is, $l_m(x) = l_m(y)$.
\item There exist points $z_1, z_p$ in $k$ such that $l_m(x) = l_k(z_1), l_m(y) = l_k(z_p),$ and $z_1 \sim z_p$. For this to be true, there must exist a sequence $z_1, z_2, \dots, z_p$ of elements of $k$ such that $l_m(x) = l_k(z_1), l_k'(z_1) = l_k'(z_2), \dots, l_k(z_p) = l_m(y)$.
\end{enumerate}

It should be clear that this set of conditions is in some way equivalent to the formula for $\theta_{A_+A_-}$ given above. By listing the conditions required for the other three possibilities (matching of points $m_+ \rightarrow n_+$, $n_- \rightarrow m_-,$ and $n_- \rightarrow n_+$) it is easy to see that they are also equivalent to the other formulas that make up the formula for $\theta$. Furthermore, using the interpretation of components of the pushout with empty preimage as loops, we observe that the loops formed in the composition correspond to equivalence classes of $\sim$ in $k$ of points that do not map to components of $A$ or $B$ whose preimages under $l_m$ or $l_n$ contain points of $m$ or $n$. This is analogous to Abramsky's definition, in which the loops formed in the composition are equivalence classes of points in $B$ that are fixed by permutations. It should then be that our definition of $\bf{1-Cob^+}$ satisfies all the properties satisfied by the definition of \cite{abram} of the free closed compact category on one element. Therefore, the composition is associative and we can define a functor from $\bf{1-Cob^+}$ to any closed compact category simply by specifying the image of the object $\{+\}$. 

Then to find a functor $\Phi: \bf{1-Cob^+} \rightarrow \bf{2-Cob}$, it suffices to prove that $\bf{2-Cob}$ is a closed compact category (as defined in \cite{abram}) and then specify the image of the object $\{+\}$. The most natural functor is that in which we map $\{+\}$ to a single boundary circle. This obviously gives rise to a functor in which both types of points in $\bf{1-Cob^+}$ are mapped to boundary circles, line segments joining points are mapped to connected components in $\bf{2-Cob}$ with two boundary circles and genus zero, and closed loops are mapped to a closed surface of genus zero. After proving the existence of this functor, we will deduce several of its properties in order to show that the orthogonal factorization system on $\bf{2-Cob}$ yields an orthogonal factorization system on $\bf{1-Cob}$. 

\subsection{2-Cob as a monoidal functor}
We first prove a lemma before showing that $\bf{2-Cob}$ is a monoidal functor in which the tensor product is the disjoint union of cobordisms.  

\newtheorem*{lem15}{Lemma 5.3.1}
\begin{lem15}
Given two spans $A_1 \leftarrow k_1 \rightarrow B_1$, $A_2 \leftarrow k_2 \rightarrow B_2$, the disjoint union of their pushouts is isomorphic to the pushout of the span $(A_1 \coprod A_2) \leftarrow (k_1 \coprod k_2) \rightarrow (B_1 \coprod B_2)$.
\end{lem15}

\begin{proof}
Let $C_1, C_2$ be the pushouts of $A_1 \leftarrow k_1 \rightarrow B_1$, $A_2 \leftarrow k_2 \rightarrow B_2$. Then given a cospan $(A_1 \coprod A_2) \rightarrow X \leftarrow (B_1 \coprod B_2)$, there is an obvious map $f$ making the diagram 
\begin{center}
\begin{tikzcd}
& k_1 \coprod k_2 \arrow{ld} \arrow{rd} \\
A_1 \coprod A_2 \arrow{rd} \arrow{rddd} && B_1 \coprod B_2 \arrow{ld} \arrow{lddd} \\
& C_1 \coprod C_2 \arrow[dashed]{dd}{f} \\ \\
& X
\end{tikzcd}
\end{center}
commute, given by the maps $C_1 \rightarrow X$, $C_2 \rightarrow X$ induced  by the cospans $A_1 \rightarrow X \leftarrow B_1, A_2 \rightarrow X \leftarrow B_2$. Now suppose there is a second distinct map $f': (C_1 \coprod C_2) \rightarrow X$ making the diagram commute. Without loss of generality, we may assume that there is some $x \in C_1$ such that $f(x) \neq f'(x)$. But then, restricting $f'$ to $C_1$, we obtain a commutative square
\begin{center}
\begin{tikzcd}
& k_1 \arrow{ld} \arrow{rd} \\
A_1 \arrow{rd} \arrow{rddd} && B_1 \arrow{ld} \arrow{lddd} \\
& C_1 \arrow[dashed]{dd}{f'} \\ \\
& X
\end{tikzcd}
\end{center}
contradicting the fact that the map $C_1 \rightarrow X$ is unique. Therefore there is a unique map $f: (C_1 \coprod C_2) \rightarrow X$ making the first diagram commute, so it is in fact the pushout.
\end{proof}

\newtheorem*{prop15_2}{Proposition 5.3.2}
\begin{prop15_2}
Given cobordisms $\phi_1: X_1 \rightarrow Y_1, \phi_2: X_2 \rightarrow Y_2, \psi_1: Y_1 \rightarrow Z_1, \psi_2: Y_2 \rightarrow Z_2$, we have $(\psi_1 \coprod \psi_2) \circ (\phi_1 \coprod \phi_2) = (\psi_1 \circ \phi_1) \coprod (\psi_2 \circ \phi_2)$. 
\end{prop15_2}

\begin{proof}
Let $\psi_i = \{m_i, k_i, A_i, a_i, b_i, g_i\}, \phi_i = \{k_i, n_i, B_i, c_i, d_i, h_i\}$. The diagrams for $\psi_1 \circ \phi_1, \psi_2 \circ \phi_2$ are

\begin{center}
\begin{tikzcd}
m_1 \arrow{rd}{a_1} && k_1 \arrow{ld}[swap]{b_1} \arrow{rd}{c_1} && n_1 \arrow{ld}[swap]{d_1} \\
& A_1 \arrow{d}{g_1} \arrow{rd}{\iota_{A_1}} && B_1 \arrow{d}{h_1} \arrow{ld}[swap]{\iota_{B_1}} \\
& \mathbb{N} & C_1 \arrow{d}{f_1} & \mathbb{N} \\
&&\mathbb{N}
\end{tikzcd}
\begin{tikzcd}
m_2 \arrow{rd}{a_2} && k_2 \arrow{ld}[swap]{b_2} \arrow{rd}{c_2} && n_2 \arrow{ld}[swap]{d_2} \\
& A_2 \arrow{d}{g_2} \arrow{rd}{\iota_{A_2}} && B_2 \arrow{d}{h_2} \arrow{ld}[swap]{\iota_{B_2}} \\
& \mathbb{N} & C_2 \arrow{d}{f_2} & \mathbb{N} \\
&&\mathbb{N}
\end{tikzcd}
\end{center}

Their disjoint union is
\begin{center}
\begin{tikzcd}
m_1 \coprod m_2 \arrow{rd} && n_1 \coprod n_2 \arrow{ld}\\
& C_1 \coprod C_2 \arrow{d}{f_1 \coprod f_2} \\
&\mathbb{N}
\end{tikzcd}
\end{center}
where the maps $(m_1 \coprod m_2) \rightarrow (C_1 \coprod C_2), (n_1 \coprod n_2) \rightarrow (C_1 \coprod C_2)$ are $\iota_{A_1}a_1 \coprod \iota_{A_2} a_2$, $\iota_{B_1}d_1 \coprod \iota_{B_2}d_2$, respectively. 

Consider the diagrams corresponding to $\psi_1 \coprod \psi_2$ and $\phi_1 \coprod \phi_2$:

\begin{center}
\begin{tikzcd}
m_1 \coprod m_2 \arrow{rd}{a_1 \coprod a_2} && k_1 \coprod k_2 \arrow{ld}[swap]{b_1 \coprod b_2} \\
& A_1 \coprod A_2 \arrow{d}{g_1 \coprod g_2} \\
&\mathbb{N}
\end{tikzcd}
\begin{tikzcd}
k_1 \coprod k_2 \arrow{rd}{c_1 \coprod c_2} && n_1 \coprod n_2 \arrow{ld}[swap]{d_1 \coprod d_2} \\
& B_1 \coprod B_2 \arrow{d}{h_1 \coprod h_2} \\
&\mathbb{N}
\end{tikzcd}
\end{center}

By Lemma 5.3.1, the diagram for $(\psi_1 \coprod \psi_2) \circ (\phi_1 \coprod \phi_2)$ is
\begin{center}
\begin{tikzcd}
m_1 \coprod m_2 \arrow{rd}{a_1 \coprod a_2} && k_1 \coprod k_2 \arrow{ld}[swap]{b_1 \coprod b_2} \arrow{rd}{c_1 \coprod c_2} && n_1 \coprod n_2 \arrow{ld}[swap]{d_1 \coprod d_2} \\
& A_1 \coprod A_2 \arrow{d}{g_1 \coprod g_2} \arrow{rd}{\iota_A} && B_1 \coprod B_2 \arrow{d}{h_1 \coprod h_2} \arrow{ld}[swap]{\iota_B}\\
&\mathbb{N}& C_1 \coprod C_2 \arrow{d}{f} & \mathbb{N} \\
&&\mathbb{N}
\end{tikzcd}
\end{center}
where $\iota_A = \iota_{A_1} \coprod \iota_{A_2}$, $\iota_B = \iota_{B_1} \coprod \iota_{B_2}$. Then the maps $(m_1 \coprod m_2) \rightarrow (C_1 \coprod C_2), (n_1 \coprod n_2) \rightarrow (C_1 \coprod C_2)$ are exactly as before. Now all that remains is to show that $f = f_1 \coprod f_2$. But this follows immediately from the fact that in the pushout square in the diagram above, elements of $k_i$ in the coproduct only map to $A_i$ and $B_i$, which in turn map only to elements of $C_i$, so the terms in the formula for computing the genus are exactly the same as in the compositions $\psi_i \circ \phi_i$.

\end{proof}

Now that we have shown $\bf{2-Cob}$ is a monoidal category, we may show it is in fact a closed compact category, which will be sufficient to prove the existence of a functor $\Phi: \bf{1-Cob^+} \rightarrow \bf{2-Cob}$.  

\newtheorem*{prop17}{Proposition 5.3.3}
\begin{prop17}
The category $\bf{2-Cob}$ is compact.
\end{prop17}
\begin{proof}
We must show that we can define a dual object $A^*$ for each object $A$ of $\bf{2-Cob}$ and that there exist maps $\eta: I \rightarrow A \otimes A^*, \epsilon: A^* \otimes A \rightarrow I$ such that the composition $\gamma^{-1}(id \otimes \epsilon)(\eta \otimes id)\gamma$ is the identity on $A$ (where $\gamma: A \rightarrow A \otimes I$ is the canonical isomorphism). For each object $A$ of $\bf{2-Cob}$, we define $A^* = A$. If $| A | = n$, let $\eta$ be given by the following diagram:
\begin{center}
\begin{tikzcd}
\emptyset \arrow{rd} && 2n \arrow{ld}[swap]{\alpha} \\
& n \arrow{d}{g = 0} \\
& \mathbb{N}
\end{tikzcd}
\end{center}
Where $\alpha$ is the map that mods out by $n$. We define $\epsilon$ as the opposite morphism. 

Now we check that the composite described above is the identity.
\begin{center}
\begin{tikzcd}
n \arrow{rd}{a} && n \coprod n \coprod n \arrow{ld}[swap]{b} \arrow{rd}{c} && n \arrow{ld}[swap]{d} \\
& 2n \arrow{rd}{\iota_1}\arrow{d}{g = 0} && 2n \arrow{ld}[swap]{\iota_2} \arrow{d}{g=0} \\
& \mathbb{N} & C\arrow{d}{g'} & \mathbb{N}\\
&& \mathbb{N}
\end{tikzcd}
\end{center}
where $a(x) = x$, $b(x,1) = x, b(x, i) = x + n$ for $i \in \{2, 3\}$, $c(x, 3) = x, c(x, i) = x+n$ for $i \in \{1,2\}$, and $d(x) = x$. 

We must determine the equivalence relation $\sim$ defined by the pushout. Note that the maps from $n \coprod n \coprod n$ to $2n$ are surjective and that by construction, an element of $n \coprod n \coprod n$ maps to an element of $2n$ equal to $x$ mod $n$ if and only if that element is $(x, i)$ for some $i$. This means that the equivalence relation on $C = A \cup B/ \sim$ identifies precisely those elements that are equal to each other mod $n$, so $C = n$. By inspection, we see that the composites map $x \in n$ to $x \in C$ along both sides of the diagram, so if we can show that the genus is zero we will have shown that this is the identity morphism on $A$. The formula for the genus from Definition 3.4 tells us that for any $x \in C$, $g'(x) = 1 - 4 + 3 = 0$. 
\end{proof}

\subsection{The functor $\bf{1-Cob^+} \rightarrow \bf{2-Cob}$}
The proof that an orthogonal factorization system on $\bf{2-Cob}$ gives an orthogonal factorization system on $\bf{1-Cob^+}$ requires a deeper understanding of the properties of the functor $\Phi: \bf{1-Cob^+} \rightarrow \bf{2-Cob}$. An orthogonal factorization system requires at least a unique factorization of every cobordism as a cofibration followed by a fibration, so a lack of faithfulness of $\Phi$ could result in a loss of this uniqueness and thus the loss of an orthogonal factorization system. Fortunately, it is in fact the case that $\Phi$ is faithful.

\newtheorem*{lem18}{Lemma 5.4.1}
\begin{lem18}
The functor $\bf{1-Cob^+} \rightarrow \bf{2-Cob}$ is faithful.
\end{lem18}
\begin{proof}
Consider the following two oriented 1-cobordisms from $m$ to $n$:
\begin{center}
\begin{tikzcd}
m \arrow{rd}{l_m} && n \arrow{ld}[swap]{l_n} \\
&C
\end{tikzcd}
\begin{tikzcd}
m \arrow{rd}{l_m'} && n \arrow{ld}[swap]{l_n'} \\
&C
\end{tikzcd}
\end{center}

In $\bf{1-Cob^+}$, we define two cobordisms as above to be equivalent if there are isomorphisms $\phi_+: m_+ \rightarrow m_+, \phi_-: m_- \rightarrow m_-, \psi_+: n_+ \rightarrow n_+, \psi_-: n_- \rightarrow n_-,$ and $\chi: C \rightarrow C$ such that 
\begin{itemize}
\item $\chi^{-1} l_m' \phi_\pm (x) = l_m(x)$
\item $\chi^{-1} l_n' \psi_\pm(x) = l_n(x)$
\end{itemize}

Now consider two 2-cobordisms between $m$ and $n$:
\begin{center}
\begin{tikzcd}
m \arrow{rd}{l_m} && n \arrow{ld}[swap]{l_n} \\
&C \arrow{d}{g} \\
& \mathbb{N}
\end{tikzcd}
\begin{tikzcd}
m \arrow{rd}{l_m'} && n \arrow{ld}[swap]{l_n'} \\
&C \arrow{d}{g'} \\
& \mathbb{N}
\end{tikzcd}
\end{center}

Two 2-cobordisms are defined to be equivalent if there exist isomorphisms $\phi: m \rightarrow m, \psi: n \rightarrow n,$ and $\chi: C \rightarrow C$ such that 
\begin{itemize}
\item $\chi^{-1} l_m' \phi(x) = l_m(x)$
\item $\chi^{-1} l_n' \psi(x) = l_n(x)$
\item $g' \chi(x) = g(x)$
\end{itemize}

We must show that two 2-cobordisms that are in the image of $\bf{1-Cob^+}$ are equivalent only if their preimages in $\bf{1-Cob}^+$ are equivalent. Note that given this definition of equivalence, the compositions of equivalent cobordisms with another cobordism yield equivalent compositions.

Suppose that we have isomorphisms $\phi, \psi, \chi$ as above, and that $\phi$ maps $s$ points of $m_+$ to $m_-$ and $t$ points of $n_+$ to $n_-$. By an elementary counting argument, there must also be $s$ points of $m_-$ mapped to $m_+$ and $t$ points of $n_-$ mapped to $n_+$. Furthermore, by the conditions on oriented 1-cobordisms, each point that is not mapped to a point with the same signature is paired with another point that is mapped to a point of different signature. This is because, for example, for any $x \in m_+$ such that $l_m(x) = c \in C$, there is exactly one point of $m_-$ that maps to $c$ or one point of $n_+$ that maps to $c$, so under $\phi$ or $\psi$ this other point must also map to a point of the opposite signature. 

For $x \in m_\pm$, there are then three cases we must consider:
\begin{enumerate}
\item $\phi(x) \in m_\pm, l_m^{-1}(l_m(x)) = \{x, x'\} \in m_\pm$
\item $\phi(x) \in m_\mp, l_n^{-1}(l_m(x)) = x' \in n_\pm$
\item $\phi(x) \in m_\mp, l_m^{-1}(l_m(x)) = \{x, x'\}, x' \in m_\mp$
\end{enumerate}

In each case, it is clear that the given $x'$ is uniquely defined. There are three analogous cases for $x \in n_\pm$, and from these we get that $x'' = x$ (by the properties of oriented 1-cobordisms). Taking these cases into consideration, we must define new isomorphisms $\phi'$ and $\psi'$ that map points to points with the same signature. If there are $p$ points $x_i \in m_+$ of type (3), then since each corresponds to a unique $x_i' \in m_-$, there are $p$ points of $m_-$ of type (3). This implies that there are precisely $s - p$ points of type (2) in both $m_+$ and $m_-$ (since the $s$ points mapped to points with different signs are necessarily of type (2) or (3)). Similarly, if there are $q$ points of type (3) in $n_+$, there are $t - q$ points of type (2) in $n_+$ and $n_-$. 

Given that each point of type (2) in $m_\pm$ corresponds to precisely one point of type (2) in $n_\pm$, it is apparent that $s - p = t - q$. Then we can label the points of type (2) in $m_+, m_-$ by $x_i, y_i$, $i = 1, \dots, s - p$, which gives us points $x'_i, y'_i$ $n_+, n_-$, respectively. Denoting the sets of points of type ($i$) by $T_i$, we can then define $\phi'$ as follows: 
\begin{center}
 $$
 \phi'(z) = \left\{
     \begin{array}{lrr}
        \phi(z) &&  z \in T_1\\
        \phi(y_i) &&  z \in T_2, \hspace{2 pt} z = x_i \\
        \phi(x_i) &&  z \in T_2, \hspace{2 pt} z = y_i \\
        \phi(z') &&  z \in T_3
     \end{array}
   \right.
$$
\end{center}

We define $\psi'$ analogously. By construction, $\phi'$ and $\psi'$ are isomorphisms that map points to points with the same signature. Furthermore, it is easy to see that any two points that mapped to the same component of $C$ are sent to a pair of points that map to the same component of $C$. We now just need to define a new isomorphism $\chi'$. For any $c \in C$ such that there is $x \in l_m^{-1}(c)$, define $\chi'(c) = l_m' \phi'(x)$. Note that this is well-defined because of the property mentioned before, that any points that map to the same component under $l_m$ to points that map to the same component under $l_m'$. Also, by construction, $l_m(x) = \chi'^{-1} l_m' \phi'(x)$. If we define $\chi'$ analogously for $c \in C$ with nonempty preimage in $n$ and let $\chi'$ be the identity for closed surfaces, we obtain the desired result. Therefore the functor $\bf{1-Cob^+} \rightarrow \bf{2-Cob}$ is faithful.
\end{proof}

To prove useful results about factorization of oriented 1-cobordisms, it will also be desirable to understand the properties of the set of intermediate boundary circles present in the cofibration-fibration factorization on $\bf{2-cob}$ and extend this to a property of factorizations of cobordisms in the image of $\Phi$. 

\newtheorem*{lem19}{Lemma 5.4.2}
\begin{lem19}
Consider a factorization of a 2-cobordism $m \rightarrow n$ of the following form:
\begin{center}
\begin{tikzcd}
m \arrow{rd}{l_m} && k \arrow{ld}[swap]{l_k} \arrow{rd}{l_k'} && n \arrow{ld}[swap]{l_n} \\
& A \arrow{rd}{\iota_A} \arrow{d}{g_1} && B \arrow{ld}[swap]{\iota_B}\arrow{d}{g_2} \\
& \mathbb{N} & C \arrow{d}{g} & \mathbb{N} \\
&& \mathbb{N}
\end{tikzcd}
\end{center}
Among all such factorizations, the cardinality of $k$ is minimal for the cofibration-fibration factorization.
\end{lem19}

\begin{proof}
Let $L_m = \iota_A l_m, L_n = \iota_B l_n$. In any factorization of this cobordism as above, by the definition of the pushout, for any element $c \in C$ in $im(L_m) \cap im(L_n)$, we must have an element $x \in k$ mapping to $c$ through $\iota_A$ and an element $y \in k$ mapping to c through $\iota_B$ (with $x, y$ not necessarily distinct). This implies that in fact $| k | \geq | im(L_m) \cap im(L_n) |$. Recall that in the cofibration-fibration factorization, $| k | = | im(L_m) \cap im(L_n) |$, so indeed $| k |$ is minimal for the cofibration-fibration factorization. 
\end{proof}

\newtheorem*{lem20}{Lemma 5.4.3}
\begin{lem20}
Suppose $\phi^+ = \{m, n, l_m, l_n, C\}$ is an oriented 1-cobordism. Then if $\phi^+ = \chi^+\psi^+, \psi^+: m \rightarrow k, \chi^+: k \rightarrow n$, $| k | \geq | im(l_m) \cap im(l_n) |$. 
\end{lem20}
\begin{proof}
The cobordisms $\chi^+, \psi^+$ map to a factorization of the image of $\phi^+$. Then it is immediate from Lemma 5.4.2 that $| k | \geq | im(l_m) \cap im(l_n) |$. 
\end{proof}

\newtheorem*{lem22}{Lemma 5.4.4}
\begin{lem22}
Let $\phi^+ = \chi^+\psi^+$ be a factorization of an oriented 1-cobordism with $\chi^+: k \rightarrow n, \phi^+: m \rightarrow k$. Then if $| k |$ is minimal among all factorizations, the isomorphism class of $k$ is uniquely determined.
\end{lem22}
\begin{proof}
Consider the following factorization of an oriented 1-cobordism:
\begin{center}
\begin{tikzcd}
m \arrow{rd}{l_m} && k \arrow{ld}[swap]{l_k} \arrow{rd}{l_k'} && n \arrow{ld}[swap]{l_n} \\
& A \arrow{rd}{\iota_A} && B \arrow{ld}[swap]{\iota_B} \\
&& C
\end{tikzcd}
\end{center}
If $\phi^+ = \{m, n, l_m, l_n, C\}$ and $| k |$ is minimal, using the factorization on $\bf{2-Cob}$ we obtain $| k | \leq | im(l_m) \cap im(l_n) |$, but Lemma 5.4.3 implies $| k | \geq | im(l_m) \cap im(l_n) |$. That is, for each pair $x \in m, y \in n$ (of the same sign) that map to the same element of $C$, there is precisely one element $z \in k$ such that $\iota_A l_k(z) = \iota_A l_m(x), \iota_B l_k'(z) = \iota_B l_n(y)$ (otherwise, we could not force $\iota_A l_m(x) = \iota_B l_n(y)$ for all pairs $x, y$ mapping to the same element of $C$). Then it must be that $k$ is a set of cardinality $| im(l_m) \cap im(l_n) |$ such that the sign of each element of $k$ matches the sign of its corresponding pair $x, y$.
\end{proof}

We will need two additional results unrelated to the lemmas above in order to prove the existence of an orthogonal factorization system on $\bf{1-Cob^+}$. 

\newtheorem*{Lem21}{Lemma 5.4.5}
\begin{Lem21}
Suppose $\phi = \chi\psi$ is a composition of cofibrations in $\bf{2-Cob}$ and that furthermore, $\phi$ and $\psi$ are each in the image of the functor $\Phi$. Then $\chi$ is the image of a cofibration in $\bf{1-Cob^+}$. Likewise, if $\phi = \chi\psi$ is a composition of fibrations with both $\phi$ and $\chi$ in the image of $\Phi$, then $\psi$ is in the image of $\Phi$.  
\end{Lem21}
\begin{proof}
We are given a composition of the following form:
\begin{center}
\begin{tikzcd}
m \arrow{rd}{a} && k \arrow[hook]{ld}[swap]{b} \arrow{rd}{c} \arrow[dashed]{ll}{\exists u} && n \arrow[hook]{ld}[swap]{d} \arrow[dashed]{ll}{\exists v}\\
& A \arrow{rd}{\iota_A} \arrow{d}{g_1} && B \arrow[hook]{ld}[swap]{\iota_B} \arrow{d}{g}\\
& \mathbb{N} & C \arrow{d}{g} & \mathbb{N} \\
&& \mathbb{N}
\end{tikzcd}
\end{center}
We must show that the functions $c, d,$ and $g_2$ satisfy the required conditions. Obviously $g_2 = 0$, since otherwise there would be an element $z \in C$ such that $g(z) \neq 0$, but we are given that the composition is in the image of $\bf{1-Cob^+}$. 

To check the required conditions for $c$ and $d$, it suffices to check that they are satisfied for each element of $B$. There are three cases to consider: for $q \in B$, either $q \in im(d)$, $q \not\in im(d)$ and $\iota_B(q) \in im(\iota_A)$, or $q \not\in im(d)$ and $\iota_B(q) \notin im(\iota_A)$. In the first case, we have $p \in C$ such that $\iota_B(q) = p$ and there exist unique elements $x \in m, y \in n$ with the same sign such that $\iota_a a(x) = \iota_B d(y) = p$. Furthermore, there is a unique element $z \in k$ with the same sign as $x$ such that $a(x) = b(z)$, and there must be some element $z' \in k$ such that $c(z') = d(y)$. 

Now we can use the fact that $g(p) = 0$. In this case, the formula for the genus simplifies to 
\begin{center}
\begin{equation*}
g(x) = 1 \hspace{8 pt} - \sum\limits_{\substack{s \in A \\ \iota_A(s) = x}} 1 \hspace{8 pt} - \sum\limits_{\substack{t \in B \\ \iota_B(t) = x}} 1 \hspace{4 pt} + \sum\limits_{\substack{y \in k \\ y \mapsto x}} 1
\end{equation*}
\end{center}
We have shown that there must be exactly one $s$ such that $\iota_A(s) = p$ and one $t$ such that $\iota_B(t) = p$, so there must be exactly one element of $k$ mapping to $p$. We already have one such element $z \in k$, however, and it has the same sign as $x$ and $y$, so we conclude $z' = z$ is the only element of $k$ mapping to $q$. 

Suppose that $q \not\in im(d)$ but $\iota_B(q) \in im(\iota_A)$. Then if $p = \iota_B(q)$, there must be two elements $x_+, x_-$ in $m_+, m_-$, respectively, such that $\iota_A a(x_\pm) = p$. Since $q$ maps to the same element of $C$ as elements of $m$, there must be some element of $k$ mapping to $p$, which means that there are $y_+, y_-$ in $k_+, k_-$ such that $a(x_\pm) = b(y_\pm)$. There can be no other element $y' \in k$ mapping to $p$, since if $\iota_A b(y') = p$, then given that $b$ is injective and the first cobordism is in the image of $\bf{1-Cob^+}$, there must be some element $x' \in m$ mapping to $p$ with $x' \neq x_\pm$. This is impossible, however, because the composition is also in the image of $\bf{1-Cob^+}$, and we already have two elements of $m$ that map to $p$. Therefore, since $\iota_B$ is injective, $y_+$ and $y_-$ must map to $q$ under the function $c$, and furthermore they must be the only elements of $k$ mapping to $q$. 

In the third case, there is nothing to check. There clearly can be no element of $k$ mapping to $q$, so $q$ is just an element of $B$ in the image of neither $c$ nor $d$. 

We have then proved that any element of $B$ satisfies the conditions imposed on an oriented 1-cobordism, so the second cobordism is in fact in the image of $\bf{1-Cob^+}$. It is clear that the analogous statement can be proved for fibrations in the exact same way. 
\end{proof}

\newtheorem*{prop3.2}{Proposition 5.4.6}
\begin{prop3.2}
The image under $\Phi$ of any cobordism in $\bf{1-Cob^+}$ factors in $\bf{2-Cob}$ as the product of two cobordisms that are also in the image of $\bf{1-Cob^+}$. 
\end{prop3.2}
\begin{proof}
Following the proof of Proposition 4.3, the diagram for an oriented 1-cobordism factors as follows:
\begin{center}
\begin{tikzcd}
m \arrow{rd}{l_m} \arrow[dashleftarrow]{rr}{\exists u} && k \arrow[hook]{ld}[swap]{b} \arrow[hook]{rd}{c} && n \arrow[two heads]{ld}[swap]{l_n} \\
& im(l_m) \sqcup (C - im(l_n) - im(l_m)) \arrow[hook]{rd}{\iota_1} \arrow{d}{g_1} && im(l_n) \arrow[hook]{ld}[swap]{\iota_2} \arrow{d}{g_2} \\
&\mathbb{N} & C \arrow{d}{g} & \mathbb{N} \\
&& \mathbb{N}
\end{tikzcd}
\end{center}
where the functions $g_1, g_2$ are zero and $k$ is the pullback of $\iota_1, \iota_2$. 

We must check that we can assign a sign to each element of $k$ such that each element of $A =  im(l_m) \sqcup (C - im(l_n) - im(l_m)), B = im(l_n)$ is mapped to by exactly two elements of the appropriate sign or no elements. Since $k$ is the pullback, as a set it is equal to $im(l_m) \cap im(l_n)$. There are then four possibilities for an element of $c \in A \cup_k B = C$: $c \in k; c \in im(l_m), c \notin k; c \in B, c \notin k, c \in A, c \notin im(l_m)$.

In the first case, $c \in A \cap B$. Since $k = im(l_m) \cap im(l_n)$, we can assign to each element $c$ the same sign as that of the elements $l_m^{-1}(c), l_n^{-1}(c)$ (which are well-defined expressions because each element of $im(l_m) \cap im(l_n)$ is mapped to by exactly one element of $m$ and $n$). The signs of $l_m^{-1}(c)$ and $l_n^{-1}(c)$ are the same by the assumption that the original cobordism is in the image of $\bf{1-Cob^+}$. After assigning signs in this way, each element $\iota_i^{-1}(c)$ (where $i \in \{1, 2\}$) is mapped to by exactly two elements of the same sign from different sets, since no additional element of $k$ can map to $\iota_i^{-1}(c)$ (given that $b, c$ are injections) and no other element of $m$ or $n$ can map to it, since each $\iota_i^{-1}(c)$ is mapped to by a unique element of $m$ or $n$ by assumption. 

In the second case, $c \notin B$, so it follows that two elements of $m$ with differing signs map to $c$, meaning that both must map to $\iota_1^{-1}(c)$. There is no element of $k$ mapping to $\iota_1^{-1}(c)$ by assumption, and there is no other element of $m$ mapping to it because then three elements of $m$ would map to a single element of $c$. 

The third case is analogous to the second.

In the fourth case, no element of $m$ or $n$ maps to $c$, so no element of $m$ maps to $\iota_1^{-1}(c)$. By assumption no element of $k$ maps to $\iota_1^{-1}(c)$, so $c$ is not in the image of $l_m$ or $b$, as desired.
\end{proof}

\newtheorem*{Thm23}{Theorem 5.4.7}

\begin{Thm23}
The orthogonal factorization system on $\bf{2-Cob}$ yields an orthogonal factorization system on $\bf{1-Cob^+}$.
\end{Thm23}
\begin{proof}
Suppose we have a square of morphisms in $\bf{1-Cob^+}$ of the following form:
\begin{center}
\begin{tikzcd}
X \arrow{r}{u} \arrow{d}{e} & X' \arrow{d}{m} \\
Y \arrow{r}{v} & Y'
\end{tikzcd}
\end{center}
where $e$ and $m$ are a cofibration and a fibration, respectively, and $u$ and $v$ are arbitrary cobordisms. This square maps to a square of the same form in $\bf{2-Cob}$. The orthogonal factorization system yields a unique 2-cobordism $w: Y \rightarrow X'$ making the diagram commute. In the proof of Theorem 4.5, $w$ is constructed by first factoring $ve = mu = fg$, where $f$ and $g$ are a fibration and a cofibration, respectively, and then noting that $g = e'e$ and $f = mm'$ for a cofibration $e'$ and a fibration $m'$. We then find $w = m'e'$. By Proposition 5.4.6, since $fg$ is in the image of $\bf{1-Cob^+}$, so are $f$ and $g$. But then Lemma 5.4.5 implies that $e'$ and $m'$ must be in the image of $\bf{1-Cob^+}$. The factorization $ve = mu = fg$ yields the diagram
\begin{center}
\begin{tikzcd}
m \arrow{rd}{a} \arrow[dashleftarrow]{rr}{\exists u} && im(l_m) \cap im(l_n) \arrow[hook]{ld}[swap]{b} \arrow[hook]{rd}{c} && n \arrow[two heads]{ld}[swap]{d} \\
& im(l_m) \sqcup (C - im(l_n) - im(l_m)) \arrow[hook]{rd}{\iota_1} && im(l_n) \arrow[hook]{ld}[swap]{\iota_2} \\
&& C
\end{tikzcd}
\end{center}
Then by Lemma 5.4.4, the isomorphism class of $k$ is uniquely determined, as $| k |$ is minimal by Lemma 5.4.3. Given that the functor $\bf{1-Cob^+} \rightarrow \bf{2-Cob}$ is faithful, the cobordisms $m \rightarrow k$ and $k \rightarrow n$ are uniquely determined, so there are unique cobordisms $f'$ and $g'$ in $\bf{1-Cob^+}$ mapping to $f$ and $g$. Since that $g = e'e$, we can write any cobordism in $\bf{1-Cob^+}$ mapping to $g$ diagrammatically as
\begin{center}
\begin{tikzcd}
m \arrow{rd}{a'} && l \arrow[hook]{ld}[swap]{b'} \arrow{rd}{c'} && k \arrow[hook]{ld}[swap]{d'} \\
& A \arrow{rd}{\iota_A} && B \arrow{ld}[swap]{\iota_B}\\
&& C
\end{tikzcd}
\end{center}
The isomorphism classes of $m$ and $l$ are known because $e$ is given, and we uniquely determined the isomorphism class of $k$ above. Then again using the fact that the functor $\bf{1-Cob^+} \rightarrow \bf{2-Cob}$ is faithful, there must be a unique cobordism from $l$ to $k$ in $\bf{1-Cob^+}$ mapping to $e'$. Similarly, the cobordism mapping to $m'$ is uniquely determined. Therefore, since $w = m'e'$, there is a unique cobordism in $\bf{1-Cob^+}$ making the original diagram commute.
\end{proof}

\section*{Acknowledgements}
I would like to thank Dr. David Spivak for introducing me to this problem and for generously helping me to solve it.

\end{document}